\documentclass[12pt]{article}
\usepackage{amsmath}
\usepackage{amsfonts}
\usepackage{amssymb}
\newcommand{\CC}{{\mathbb C}}
\newcommand{\RR}{{\mathbb R}}
\newcommand{\HH}{{\mathbb H}}
\newcommand{\BB}{{\mathbb B}}
\newcommand{\PP}{{\mathbb P}}

\newcommand{\ra}{\rightarrow}
\newcommand{\ov}{\overline}
\newcommand{\p}{\partial}

                      \newtheorem{theorem}{Theorem}[section]
\newtheorem{lemma}[theorem]{Lemma}

\date{October 25, 2007}

\begin{document}

\author{ {\it Dedicated to Professor  Kohn on the occasion of his 75th
birthday}\\ \ \\
 Xiaojun  Huang\footnote{
Supported in part by NSF-0500626},\  Shanyu Ji and Yuan Zhang}
\title{\bf A criterion  for a proper rational  map to be equivalent to a proper polynomial  map}
 \maketitle
\section{Introduction}
Let $\BB^{n}$ be the unit ball in the complex space $\CC^{n}$.
Write $Rat(\BB^{n}, \BB^{N})$ for the space of  proper rational
holomorphic maps from $\BB^{n}$ into $\BB^{N}$
and $Poly(\BB^{n}, \BB^{N})$ for the set of proper polynomial
holomorphic maps from $\BB^{n}$ into $\BB^{N}$.
 We say that $F$ and $G\in Rat(\BB^{n}, \BB^{N})$ are
{\it equivalent} if there are automorphisms $\sigma\in Aut(\BB^{n})$
and $\tau\in Aut(\BB^{N})$ such that $F=\tau\circ G\circ\sigma$.

Proper rational holomorphic maps from ${\BB}^n$ into ${\BB}^{N}$
with $N\le 2n-2$ are equivalent to the identity map ([Fa] [Hu]).
In [HJX], it is shown that $F\in Rat({\BB}^n,{\BB}^N)$ with $N\le
3n-4$ is equivalent to a quadratic monomial map, called the
D'Angelo map. However, when the codimension is sufficiently large,
there are a plenty of rooms to construct rational holomorphic maps
with certain arbitrariness by the work in Catlin-D'Angelo  [CD].
Hence, it is reasonable to believe that after lifting the
codimension restriction,  many proper rational holomorphic maps are not
equivalent to polynomial proper holomorphic maps. In  the last paragraph
of the paper [DA], D'Angelo gave a philosophic discussion on this
matter. However, explicit examples of  proper rational holomorphic
maps, that are not equivalent to  polynomial proper holomorphic
maps, do not seem to exist
 in the literature.
And the problem of  determining if an explicit proper rational
holomorphic map is equivalent to a polynomial holomorphic map
 does not seem to have been  studied so far.

This short paper  is concerned with such a problem. We will first
give an explicit criterion when a rational holomorphic map is
equivalent to a polynomial holomorphic map. Making use of this
criterion, we construct rational holomorphic maps of degree 3 that
are `almost' linear but are not equivalent to polynomial
holomorphic maps. On the other hand, with the help of the
classification result  in \cite{CJX}, our criterion is used in
this paper  to show that any proper rational holomorphic map from
${\BB}^2$ into ${\BB }^N$ of degree two is equivalent to a
polynomial holomorphic map.

\medskip
{\bf Acknowledgment}: Part of our work was done when the authors
were participating the AIM workshop on Complexity of Mappings in CR
Geometry in 2006. The authors would like to thank J. D'Angelo and P.
Ebenfelt for the invitation and the AIM for having provided to them
a   stimulating working environment. The second author would  like
to thank J. Faran for several stimulating discussions related to the
topics of this paper. The third author would  also like to thank F.
Meylan for some discussions related to the paper.

\section {A criterion}
 Let $F=\frac{P}{q}=\frac{(P_1, ..., P_N)}{q}$ be
a rational holomorphic map from the unit ball $\BB^n \subset \CC^n$
into the unit ball  $\BB^N \subset \CC^N$, where $(P_j)_{j=1}^N$,
$q$ are polynomial holomorphic functions and $(P_1, ..., P_N, q)=1$.
We define $\deg(F)=\max\{deg(P_j)_{j=1}^N, deg(q)\}$. Then $F$
induces a rational  map from $\CC\PP^n$ into $\CC\PP^N$ given by
\[
\hat F([z_1:...:z_n:t]) = \bigg[t^k P(\frac{z}{t}): t^k q(\frac{z}{t})
\bigg]
\]
where $z=(z_1, ..., z_n)\in \CC^n$ and $\deg(F)=k$.

$\hat F$ may not be holomorphic in general. For instance, we have
the following:

\bigskip

\noindent{\bf Example 2.1} \\ I. Let  $F_\theta(z, w)=(z, cos
\theta\ w, sin \theta\ zw, sin \theta\ w^2)$ be the proper
monomial map from ${\BB}^2$ into ${\BB}^4$ (called the D'Angelo
map), where $0<\theta<\frac{\pi}{2}$. Then the pole of
$\hat{F_\theta}$ consists of one point: $\{ [z: w: t]\ |\ w=0,\
t=0\} = \{[1:0:0]\}$.

\noindent II. Let $G_\alpha = (z^2, \sqrt{1+cos^2 \alpha}\ zw,
cos\alpha\ w^2,sin\alpha\ w)$ be the proper monomial map from
${\BB}^2$ into ${\BB}^4$ where $0\le \alpha <\frac{\pi}{2}$. Then
$G_\alpha$ induces
\[
\hat G_\alpha =  [ z^2: \sqrt{1+cos^2 \alpha}\ zw: cos\alpha\
w^2:sin\alpha\ wt: t^2].\] There is no pole for $\hat G_\alpha$.
Hence $\hat G_\alpha$ is holomorphic.

\bigskip

Write $\BB^n_1 = \{[z_1: ...:z_n : t]\in \CC\PP^n\ |\
\sum^n_{j=1}|z_j|^2<|t|^2 \}$, which is the projectivization of
$\BB^n$. Write $U(n+1, 1)$ for the collection of the linear
transforms $A$ such that
\[
A E_{n+1, 1} \ov A^t = E_{n+1, 1}\]
where
\[
E_{n+1, 1} = \begin{pmatrix}
I_n & 0 \\
0 & -1
\end{pmatrix}.
\]
Then $U(n+1, 1)/ \{\pm Id\} = Aut(\BB^n_1)\approx Aut(\BB^n)$.

\begin{lemma}
\label{lemma}
For any hyperplane $H\subset \CC\PP^n$ with $H\cap \ov{\BB^n_1}=\emptyset$,
there is a
$\sigma \in U(n+1, 1)$ such that $\sigma(H)=H_\infty = \{[z_1:...:z_n: 0]\}$.
\end{lemma}

\noindent{\it Proof:}\ \ Assume that $H: \sum^n_{j=1} a_j z_j -
a_{n+1} t = 0$ with $\vec a =(a_1, ..., a_{n+1})\not=0$. Under the
assumption that $H\cap \ov{\BB^n_1}=\emptyset$, we have
$a_{n+1}\not=0$. Without loss of generality, we can assume
$a_{n+1}=1$. Let $U$ be an $n\times n$ unitary matrix such that
\[
(a_1, ..., a_n) U = (\lambda, 0, ..., 0).\] Let $\sigma =
\begin{pmatrix}
U & 0 \\
0 & I\\
\end{pmatrix}$.
Then $\sigma(H)=\{[z: t] \in \CC\PP^n\ |\ \lambda z_1 - t=0\}$
with $|\lambda|<1$. Let $T\in Aut(\BB^n)$ be defined by
\[
T(z,z') = \bigg(\frac{z_1 -\ov\lambda }{1- \lambda z_1},\
\frac{\sqrt{1-|\lambda|^2} z'}{1-\lambda z_1} \bigg)\] with
$z'=(z_2, ..., z_n)$. Then $\hat T\in U(n+1, 1)$ is defined by
\[
\hat T([z: t]) = [z_1 - \ov \lambda t: \sqrt{1 - |\lambda|^2} z':
t - \lambda z_1].\] Then $\hat T\circ \sigma$ maps $H$ to
$H_\infty$. \ $\Box$

Our criterion can be stated as follows:

\begin{theorem}
\label{main them} Let $F\in Rat(\BB^n, \BB^N)$. Then $F$ is
equivalent to a proper polynomial holomorphic map from $\BB^n$ to
$\BB^N$ if and only if there exist hyperplanes $H\subset \CC\PP^n$
and $H'\subset \CC\PP^N$ such that $H\cap \ov{\BB^n_1}=\emptyset$,
$H'\cap \ov{\BB^N_1} = \emptyset$ and
\[
\hat F(H) \subset H',\ \ \hat F(\CC\PP^n \backslash  H) \subset
\CC\PP^N \backslash H'.\]
\end{theorem}

\noindent{\it Proof:}\ \ If $F$ is a polynomial holomorphic map,
then $\hat F=[t^k F(\frac {z}{t}), t^k]$ with $deg(F)=k$. Let
$H=H_\infty$ and $H'=H_\infty'$. Then they satisfy the property
defined in the lemma.

If $F$ is equivalent to a polynomial holomorphic map $G$, then
there exist $\hat \sigma\in U(n+1, 1), \hat \tau\in U(n+1, 1)$
such that $\hat F=\hat \tau \circ \hat G \circ \hat \sigma$. Let
$H=\hat \sigma^{-1}(H_\infty)$ and $H'=\hat\tau(H_\infty')$. Then
they are the desired ones.

Conversely, suppose that $\hat F,$ $H$ and $H'$ are as in the
theorem. By Lemma \ref{lemma},  we can find  $\hat \sigma\in
U(n+1, 1)$ and $\hat \tau\in U(n+1, 1)$ such that $\hat
\sigma(H)=H_\infty$  and $\hat \tau(H')=H'_\infty$. Let $\hat
Q=\hat \tau\circ \hat F\circ \hat \sigma^{-1}$. Then $\hat Q$
induces a proper rational holomorphic map $Q$ from $\BB^n$ into
$\BB^N$. If $Q=\frac{P}{q}$ where $(P, q)=1$ and $deg(Q) =k$, then
\[
\hat Q = [t^k P(\frac{z}{t}) : t^k q(\frac{z}{t})].
\]

Suppose  that $q\not\equiv constant$. Let $z_0\in \CC^n$ be such
that $q(z_0)=0$. Notice that $\hat Q(H_\infty)\subset H'_\infty$ and
$\hat Q(\CC\PP^n \backslash H_\infty) \subset \CC\PP^N\backslash
H'_\infty$. However $t^k q(\frac{z_0}{t}) = 0$ for $t=1 $. Hence,
$\hat Q([ z_0: 1])\subset H'_\infty$. That gives  a contradiction.
$\Box$

\bigskip

Write the Cayley transformation
\[
\rho_n(z', z_n) = \bigg( \frac{2z'}{1- i z_n},\ \frac{1+i z_n}{1-
iz_n} \bigg).\] Then $\rho_n$ biholomorphically maps $\p\HH^n$ to
$\p\BB^n \backslash \{(0, 1) \}$. $\rho_n$ induces a linear
transformation of $\CC\PP^n$:
\[
\hat \rho_n = [2 z' : t+i z_n : t- i z_n].\]
$\hat \rho_n$ maps ${\bf S}^n_1=\{[z: t] \in \CC\PP^n\ |\ \frac{z_n \ov t - t
\ov z_n}{2i}
> |z'|^2\}$ to $\BB^n_1$.

Now let $G$ be a CR map from $\p\HH^n$ to $\p\HH^N$. Then $\rho_N
\circ G\circ \rho_n^{-1}$ extends to a proper rational holomorphic
map from $\BB^n$ to $\BB^N$. Then, Theorem \ref{main them} can be restated as:

\begin{theorem}
$\rho_N\circ G\circ \rho_n^{-1}$ is equivalent to a polynomial
holomorphic map if and only if there are $H\subset \CC\PP^n$,
$H'\subset \CC\PP^N$ such that $\hat G(H)\subset H'$ and $\hat
G^{-1}(H')\subset H$ with
\[
H\cap \ov{{\bf S}^n_1} = \emptyset,\ \ H'\cap \ov{{\bf S}^N_1} = \emptyset. \]
\end{theorem}

\section{proper rational holomorphic maps from ${\BB}^2$ into ${\BB}^N$ of degree two}
As a first application of Theorem 2.2, we prove the following:

\begin{theorem}
A map $F\in Rat({\BB}^2,{\BB}^N)$ of degree two is equivalent to a
polynomial proper holomorphic map in $Poly({\BB}^2,{\BB}^N)$.
\end{theorem}

\noindent{\it Proof:}\ \ By \cite{HJX}, we know that any rational
holomorphic map of degree 2 from ${\BB}^2$ into ${\BB}^N$ is equivalent to a
map of the form $(G,0)$, where
the map $G$ is from ${\BB}^2$ into ${\BB}^5 $. Hence, to prove Theorem 3.1, we need only to assume that $N=5$.
 After applying a Cayley transformation and using the result in \cite{CJX},
we can assume that $F=(f, \phi_1, \phi_2$, $\phi_3, g)$ is from
${\HH}^2$ into ${\HH}^5$ with either \\(I)
\begin{eqnarray*}
 f=\frac{z + \frac{i}{2}zw}{1+e_2 w^2},\ \phi_1 =\frac{z^2}{1 +
e_2 w^2}, \ \phi_2=\frac{c_1 zw}{1 + e_2 w^2},\ \phi_3=0,\ g
=\frac{w}{1 + e_2 w^2}
\end{eqnarray*} where $-e_2=\frac{1}{4}+c_1^2$ and $c_1>0$ or\\
(II)
\[
f=\frac{z + (\frac{i}{2} + i e_1)zw}{1+i e_1 w + e_2 w^2},\
\phi_1 =\frac{z^2}{1+i e_1 w + e_2 w^2}, \]
\[
\phi_2=\frac{c_1 zw}{1+i e_1 w + e_2 w^2},\ \phi_3 =\frac{c_3
w^2}{1+i e_1 w + e_2 w^2}, \ g =\frac{w+i e_1 w^2}{1+i e_1 w + e_2
w^2}\] where $-e_1, -e_2 >0$, $c_1, c_3>0$ and
\[
e_1 e_2 = c_3^2,\ \  -e_1-e_2=\frac{1}{4} + c_1^2. \]
Write
$[z:w:t]$ for the homogeneous coordinates of $\CC\PP^2$, then in
Case (I) it induces a meromorphic map $\hat F: \CC\PP^2 \ra
\CC\PP^5$ given by
\[\hat F([z:w:t]) = [tz +\frac{i}{2}zw: z^2 : c_1 zw: 0: tw : t^2+ e_2 w^2]\ \ \ \forall [z:w:t]\in \CC\PP^2,\]
and in Case (II), it induces a meromorphic map $\hat F: \CC\PP^2
\ra \CC\PP^5$ given by
\[
\hat F([z:w:t]) = [tz +(\frac{i}{2} + i e_1)zw: z^2 : c_1 zw: c_3
w^2: tw + i e_1 w^2: t^2+i e_1 wt + e_2 w^2]\]$\forall [z:w:t]\in
\CC\PP^2$.\\ In terms of Theorem 2.3, we will look for $H=\{ - t =
\mu_1 z_1 + \mu_2 z_2 \} \subset \CC\PP^2$ and $H'=\{ - t' =
\sum^5_{j=1} \lambda_j z_j'\} \subset \CC\PP^5$ such that $H\cap
\ov{{\bf S}^2_1} = \emptyset,\ \ H'\cap \ov{{\bf S}^5_1} =
\emptyset$ with
\[
\hat F(H)\subset H'\ \  and \ \ \hat F^{-1}(H')\subset H.\] \\
We next prove the following lemma:

\begin{lemma}
Let $H=\{ -t = \sum^n_{j=1}K_j z_j  \} \subset \CC\PP^n$. Then
$H\cap \ov{{\bf S}^n_1} = \emptyset$ if and only if
\[4 \Im( K_n ) + \sum^{n-1}_{j=1} |K_j|^2 < 0.\]
\end{lemma}

\noindent{\it Proof:}\ \ \ \ Suppose for $z_j $ and $t= -
\sum^n_{j=1} K_j z_j$, we have
\[
\frac{w \ov t - t \ov w}{2i} < \sum^{n-1}_{j=1}|z_j|^2.\] Here we
identify $z_n = w$. We then get
\[
\frac{- \ov{K_n}|w|^2 + K_n |w|^2}{2i} + \sum^{n-1}_{j=1} \frac{-
\ov{K_j} \ov {z_j} w + K_j z_j \ov w}{2i} < \sum^{n-1}_{j=1}
|z_j|^2.\] Hence
\[
|w|^2 \Im(K_n) < \sum^{n-1}_{j=1} \{|z_j|^2 - 2
\Re(\frac{K_j}{2i}z_j \ov{w} ) \},\] or
\[
|w|^2\bigg(\Im(K_n) + \sum^{n-1}_{j=1} \frac{|K_j|^2}{4} \bigg) <
\sum^{n-1}_ {j=1} |z_j - \frac{\ov{K_j}}{2i}w |^2.\] Since $\{z_j,
w\}$ are independent variables, this can only happen if and only
if
\[
\Im(K_n) + \sum^{n-1}_{j=1} \frac{|K_j|^2}{4}<0.\] This proves the
lemma. $\Box$

\medskip

We remark that for the map $\hat F$ as defined above, if $H$ and
$H'$ are hyperplanes sitting in $\CC\PP^2$ and $\CC\PP^5$,
respectively, such that
\[
H \cap {\bf S}^2_1 = \emptyset,\ H'\cap {\bf S}^5_1 = \emptyset,\
\ \hat F(H)\subset H', \ \hat F^{-1}(H')\subset H,\] then $H$ and
$H'$ have to be defined, respectively,  by  equations of the form:
\[
H: -t = \mu_1 z + \mu_2 w,\ \ H':- t'=\sum^5_{j=1} \lambda_j z_j'.
\]
Indeed, if $H'=\{\sum_{j=1}^5 \lambda_j z_j'=0 \}$, then
\[
\lambda_1[tz + i(\frac{1}{2}+e_1)zw ] + \lambda_2 z^2 + \lambda_3
c_1 zw + \lambda_4 c_3 w^2 + \lambda_5(tw + i e_1 w^2) = (\mu_0 t
+ \mu_1 z + \mu_2 w) ^2\]$\forall [z:w:t]\in\CC\PP^2$.\\ We notice
that $\lambda_5\not=0$; otherwise $H'\cap {\bf S}
^5_1\not=\emptyset$. We thus get $\mu_0=0$. This would immediately
give a contradiction. Similarly, if $H=\{ \mu_1 z + \mu_2 w=0\}$,
then it follows that $H'$ is defined by $\{ \sum_{j=1}^5 \lambda_j
z_j'=0\}$. We also reach a
contradiction.\\

Summarizing the above, we have:\medskip

\textbf{Case (I)} We need only to find out $\mu_1, \mu_2,
\lambda_1, ..., \lambda_5\in \CC$ such that
\[4\Im(\mu_2)+|\mu_1|^2<0,\ \ \ \ \ \
4\Im(\lambda_5)+\sum_{j=1}^4|\lambda_j|^2<0\] and
\begin{eqnarray*} \lambda_1(tz +\frac{i}{2} zw)+\lambda_2z^2
+\lambda_3 c_1 zw+\lambda_5  tw + (t^2 + e_2
w^2)=(t+\mu_1z+\mu_2w)^2 \ \ \forall[z:w:t]\in
\CC\PP^2.\end{eqnarray*}

 It is easy to verify that  $\lambda_1=\lambda_2=\lambda_3=\lambda_4=\mu_1=0$,
$\lambda_5=-2\sqrt{|e_2|}i$ and $\mu_2=-\sqrt{|e_2|}i$ satisfy the
above conditions. Hence in  case (I), the map  is always
equivalent to
a polynomial holomorphic map.\\

\medskip
\textbf{Case (II)} Similar to case (I),  we need  to find out
$\mu_1, \mu_2, \lambda_1, ..., \lambda_5\in\CC$ such that
\[4\Im(\mu_2)+|\mu_1|^2<0,\ \ \ \ \ \
4\Im(\lambda_5)+\sum_{j=1}^4|\lambda_j|^2<0\] and
\begin{eqnarray*}
&&
\lambda_1 (tz + i(\frac{1}{2} + e_1) zw) + \lambda_2 z^2 + \lambda_3 c_1 zw +
\lambda_4 c_3 w^2 + \lambda_5(tw+ie_1 w^2)\\
&&+(t^2+i e_1 tw + e_2 w^2) \equiv (t + \mu_1 z + \mu_2 w)^2,\ \ \forall [z:
w: t]
\in \CC\PP^2.
\end{eqnarray*}

Comparing the coefficients, we get
\begin{eqnarray*}
&&\lambda_1= 2 \mu_1,\ \lambda_2=\mu_1^2,\
\lambda_3=\frac{1}{c_1}[- i (1+2e_1) \mu_1 + 2\mu_1 \mu_2],\\
&&\lambda_4=\frac{1}{c_3}(\mu_2^2 - e_2 - 2i e_1 \mu_2 - e_1^2 ),\
\lambda_5=2\mu_2 - i e_1.
\end{eqnarray*}
In sum, we obtain the following statement:

 \medskip

 {\it $\rho_N\circ F \circ \rho_n^{-1}$
is equivalent to a polynomial holomorphic map if and only if there
are $\mu_1, \mu_2 \in \CC$ such that $4 \Im(\mu_2) + |\mu_1|^2 <
0$ and that
\[
-4 e_1 + 8 \Im(\mu_2) + 4|\mu_1|^2 + |\mu_1|^4 + \frac{1}{c_1^2}
|2\mu_1 \mu_2 - i(1+2e_1) \mu_1|^2 + \frac{1}{c_3^2} |\mu_2^2 - e_2
-e_1^2 - 2i e_1 \mu_2|^2 < 0.\]}

 We will look for $\mu_1$ and $\mu_2$  with
$\mu_1 = 0$ and $\mu_2=iy$ ($y<0$).

To prove that $\rho_N\circ F \circ \rho_n^{-1}$ is equivalent to a
polynomial holomorphic map, it suffices for us to show that there
exists $y < 0$ such that
\[
-4 e_1 + 8 y + \frac{1}{c_3^2}(-y^2 - e_2 - e_1^2 + 2 e_1 y)^2 <0,
\] or
\[
J(y):=(-4 e_1 + 8y)e_1 e_2 + (y^2 - 2 e_1 y + e_1^2 + e_2)^2
= (8y - 4 e_1) e_1 e_2 + ((y-e_1)^2 + e_2)^2 < 0.\]
Notice that as a function in $y<0$,
\[
\lim_{y\ra -\infty} J(y) = + \infty,\ J(0)=(e_1^2 - e_2)^2 >0.\]
We need to show that
\[
\min_{y\le  0} J(y) < 0. \]

Notice that $J'(y)= 8 e_1 e_2 + 4((y-e_1)^2+e_2)(y-e_1)$. Setting
$J'(y)=0$, we get
\[
(y-e_1)^3 + e_2(y-e_1) + 2 e_1 e_2 = 0.\] $J'(y)=0$ thus has a root
$y_0\in (-\infty, 0)$; for
\[
\lim_{y\ra - \infty} J'(y) = -\infty,\ \ J'(0)=4(-e_1^3 + e_1 e_2)>0. \]
Let $\zeta_0, \zeta_1, \zeta_2$ be the solution of
\[ \zeta^3 + e_2 \zeta + 2 e_1 e_2 = 0\ \ with\ \zeta_0 = y_0 - e_1.\]
Then $\zeta_0 + \zeta_1 + \zeta_2 = 0$, $\zeta_0 \zeta_1 + \zeta_0
\zeta_2 + \zeta_1 \zeta_2 = e_2$ and $\zeta_0 \zeta_1 \zeta_2 = -
2 e_1 e_2$. Hence $\zeta_0 =- \zeta_1 - \zeta_2$. We get
\[
- \zeta^2_0 + \zeta_1 \zeta_2 = e_2,\] or $\zeta_1 \zeta_2 = e_2 +
\zeta_0^2$, and
\[
\frac{1}{\zeta_1 \zeta_2}=-\frac{\zeta_0}{2e_1e_2}.\] In
particular, $\frac{1}{\zeta_1 \zeta_2} \in \RR\backslash\{0\}. $

Now $J(y_0) = ( - 4 e_1 + 8 \zeta_0 + 8 e_1)e_1 e_2 + (\zeta_0^2 +
e_2)^2$ = $2 e_1 e_2 (4 \zeta_0 + 2e_1) + (\zeta_1 \zeta_2)^2$ =
$- \zeta_1 \zeta_1 \zeta_2(4 \zeta_0 + 2 e_1) + (\zeta_1
\zeta_2)^2$.

Notice that $4 \zeta_0^3 = - 8 e_1 e_2 - 4 e_2 \zeta_0$. We see
that
\begin{eqnarray*}
&& 2e_1e_2\frac{J(y_0)}{(\zeta_1 \zeta_2)^2} = 2 e_1 e_2 +
\zeta_0^2 (4 \zeta_0 + 2 \zeta_1 )
= 2 e_1 e_2 - 8 e_1 e_2 - 4 e_2 \zeta_0 + 2e_1 \zeta_0^2\\
&&
= - 6 e_1 e_2 - 4 e_2 \zeta_0 + 2 e_1 \zeta_0^2
= - 2 e_2 (3e_1 + 2 \zeta_0) + 2 e_1 \zeta_0^2.
\end{eqnarray*}
Since $\zeta_0 = y_0 - e_1 < - e_1$, $3 e_1 + 2\zeta_0 < e_1 <0$.
Therefore $\frac{J(y_0)}{(\zeta_1 \zeta_2)^2} 2 e_1 e_2 < 0$. Here
we showed that $J(y_0)<0$. This completes the proof of Theorem
3.1. $\Box$


\bigskip

Our proof of Theorem 3.1 is, in fact, a constructive proof, which
can be used to find precisely the  polynomial holomorphic maps
equivalent to the original ones. In the following, we demonstrate
this by giving an explicit example:
\bigskip

\noindent{\bf Example 3.3}\ \  Let
$F=(f,\phi_1,\phi_2,\phi_3,g):\HH^2\rightarrow \HH^5$ be defined
as follows:
\[ f(z,w)=\frac{z - \frac{i}{2} zw}{1-i  w - \frac{1}{3}
w^2},\ \phi_1(z,w) =\frac{z^2}{1-i  w -\frac{1}{3} w^2}, \]
\[
\phi_2(z,w)=\frac{\sqrt{\frac{ 13}{12}} zw}{1-i w -\frac{1}{3}
w^2},\ \phi_3(z,w) =\frac{\frac{\sqrt{3}}{3} w^2}{1-i  w
-\frac{1}{3} w^2}, \ g(z,w) =\frac{w-i  w^2}{1-i  w -\frac{1}{3}
w^2}\] It is equivalent to the proper polynomial holomorphic map
$G$ from ${\BB}^2$ into ${\BB}^5$:
\[G(z,w)=\bigg{(}\frac{\sqrt 3}{9}(-2+4z+z^2), -\frac{\sqrt
6}{9}(1+z+z^2),\frac{\sqrt 3}{12}(5+3z)w,\frac{\sqrt
6}{6}w^2,\frac{\sqrt {13}}{12}i(1-z)w\bigg{)}.
\]
\noindent{\it Proof:}\ \ \ \ In fact, for the map $F$ given above,
$e_1=-1,\ e_2=-\frac{1}{3},\ c_1=\sqrt{\frac{13}{12}},\
c_3=\frac{\sqrt 3}{3}$. From the proof of Theorem  3.1, $\hat
F(H)=H'$ where $H \subset\CC\PP^2,\ H'\subset\CC\PP^5$  are
defined by
\begin{eqnarray*}
&& H: t = - y_0 i w,\ or\ \frac{w}{t}=\frac{i}{y_0},\\
&& H': t' = - \lambda_4 z_4' - \lambda_5 w', \ or\  - \lambda_4
\frac{z_4'}{t'} - \lambda_5 \frac{w'}{t'}=1.
\end{eqnarray*}
Here $y_0<0$ is a solution for $(y_0 +1)^3 -\frac{1}{3}(y_0 +1) +
\frac{2}{3} = 0,\ \lambda_4 = \frac{1}{c_3}[-(y_0 - e_1)^2 - e_2]
=- \frac{(y_0 - e_1)^2 + e_2} {\sqrt{e_1 e_2}}$ and $\lambda_5 =
2i y_0 - e_1 i$. Therefore $y_0=-2$,  $\lambda_4=-\frac{2}{\sqrt
3}$ and $\lambda_5=-3i$. Thus we see that
\begin{eqnarray*}
&& H: t = 2iw,\ or\ \frac{w}{t}=\frac{1}{2i},\\
&& H': t' = \frac{2}{\sqrt 3} z_4' +3i w' \ or\ \frac{2}{\sqrt
3}\frac{z_4'}{t'} + \frac{3iw'}{t'}=1.
\end{eqnarray*}
  Consider
$\tilde{F}:=\rho_5\circ F\circ \rho_2^{-1}: \BB^2\rightarrow \BB^5$
where $\rho_i$ are the corresponding Cayley transformations. An easy
computation shows  that the projectivization of $\tilde F$, denoted
by $\hat{ \tilde F}$, is as follows:
 \begin{eqnarray*}
 \hat {\tilde
 F}([z:w:t])=\bigg[z(3t+w):2z^2:2i\sqrt{\frac{13}{12}}z(t-w):-\frac{2\sqrt
 3}{3}(t-w)^2\\:
\frac{1}{3}(t^2+10tw+w^2):\frac{1}{3}(13t^2-2tw+w^2)\bigg]
\end{eqnarray*}
and
\begin{eqnarray*}
&& \hat{\tilde H}:= \hat{\rho}_2(H): t = \frac{1}{3}w,\\
&& \hat{\tilde {H'}}= \hat{\rho}_5(H'): t' = \frac{\sqrt 3}{6}
z_4' +\frac{1}{2} w'.
\end{eqnarray*}
Clearly, we have $\hat{\tilde H}\subset \CC\PP^2$ and $\hat{\tilde
{H'}}\subset \CC\PP^5$ satisfying the property that $\hat{\tilde
H}\cap \ov{\BB^2_1}=\emptyset$, $\hat{\tilde {H'}}\cap \ov{\BB^5_1}
= \emptyset$ and
\[
\hat {\tilde F}(\hat {\tilde H}) \subset \hat {\tilde {H'}},\ \
\hat {\tilde F}(\CC\PP^2 \backslash \hat {\tilde H}) \subset
\CC\PP^5 \backslash \hat {\tilde {H'}}.\] According to Lemma
\ref{lemma}, let
\begin{eqnarray*}
\hat\sigma_1([z:w:t])&=&\bigg[\frac{2\sqrt
2}{3}w:z+\frac{t}{3}:t+\frac{z}{3}\bigg]\\
\hat\sigma_2([z'_1:z'_2:z'_3:z'_4:w':t'])&=&\bigg[\frac{1}{2}(z'_4+\sqrt
3 w')-\frac{\sqrt 3}{3}t:\frac{\sqrt 6}{6}(w'-\sqrt
3z'_4)\\&&:\frac{\sqrt 6}{3}z_1:\frac{\sqrt 6}{3}z_2:\frac{\sqrt
6}{3}z_3:t-\frac{\sqrt 3}{6}(z'_4+\sqrt 3 w')\bigg],
\end{eqnarray*}
then $\hat\sigma_1\in U(3,1)$ and $\hat\sigma_2\in U(6,1)$ with
$\hat\sigma_1(\hat{\tilde H}_{\infty})=\hat{\tilde H}$ and
$\hat\sigma_2(\hat{\tilde {H'}})=\hat{\tilde {H'}}_{\infty}$. The
desired proper polynomial holomorphic map  $G$  is thus given  by
$\hat\sigma_2\circ\hat {\tilde F}\circ\hat\sigma_1$. $\Box$

\medskip
{\bf  Remark 3.4}: It may be interesting  to notice that  the map
$G$ in Example 3.3 does not preserve the origin and does not equivalent to a map
of the form $(G',0)$.
 We do not know
other examples of polynomial proper  holomorphic maps between balls of this type.

\section{Non-polynomially equivalent proper rational holomorphic maps}
In this section, we apply Theorem 2.2 to construct examples of
rational holomorphic maps which are not equivalent to polynomial
holomorphic maps.

\noindent{\bf Example 4.1}: Let $G(z,w)= \bigg(z^2, \sqrt{2}zw,
w^2(\frac{z-a}{1- \ov a z},  \frac{\sqrt{1- |a|^2}w}{1- \ov a z})
\bigg)$, $|a|<1$, be a map in $Rat(\BB^2$, $\BB^4)$. Then $G$ is
equivalent to a proper polynomial holomorphic map in $Poly({\BB}^2,{\BB}^4)$ if and only if $a=0$.
\bigskip

\noindent{\it Proof:}\ \ \ \ Indeed, we have
\[\hat G=\bigg[(t-\ov a z)z^2:(t-\ov a z)\sqrt 2 zw:w^2(z-at):w^2\sqrt{1-|a|^2}w:(t^3-\ov a
t^2z)\bigg].\] Suppose there exist hyperplanes
$H=\{\mu_1z_1+\mu_2w+\mu_0t=0\}\subset \CC\PP^2$ and $H'=\{
\sum_{j=1}^4\lambda_jz'_j+\lambda_0t'=0\}\subset \CC\PP^4$ such
that $H\cap \ov{{\bf S}^2_1}=\emptyset$, $H'\cap \ov{{\bf S}^4_1}
= \emptyset$ and $ \hat F(H) \subset H',\ \ \hat F(\CC\PP^2
\backslash  H) \subset \CC\PP^4 \backslash H'$. Then
\begin{eqnarray*}
\lambda_1(t-\ov a z)z^2 +\lambda_2 (t-\ov a z)\sqrt 2 zw+\lambda_3
w^2(z-at)+ \lambda_4w^2\sqrt{1-|a|^2}w \\+\lambda_0(t^3-\ov a
t^2z) =(\mu_1z+\mu_2w+\mu_0t)^3\ \ \
\forall[z:w:t]\in\CC\PP^2.\end{eqnarray*}
 Apparently
$\lambda_0\neq0$. Hence we can assume that $\lambda_0=1, \mu_0=1$.
Therefore by comparing the coefficients we get
\begin{eqnarray*}
\mu_1^3=-\ov a\lambda_1,\ \mu_2^3=\lambda_4\sqrt{1-|a|^2},\
3\mu_2=0,\ 3\mu_1=0\\
3\mu_1^2=\lambda_1,\ 6\mu_1\mu_2=\sqrt 2\lambda_2,\
3\mu_1^2=\lambda_3=0.
\end{eqnarray*}
We then have
$\lambda_1=\lambda_2=\lambda_3=\lambda_4=\mu_1=\mu_2=0$. Moreover,
the above equality can only hold when $a=0$. By Theorem 2.2, we see the conclusion.
$\Box$

\bigskip

\noindent{\bf Example 4.2}: Let $ F(z',w)=\bigg(z', wz',
w^2(\frac{\sqrt{1-|a|^2}z'}{1- \ov a w}, \frac{w - a}{1- \ov a w})
\bigg)$ with $|a|<1$ be a map in $Rat(\BB^n, \BB^{3n-2})$. Then
$F$ has geometric rank $1$ and is linear along each hyperplane
defined by $w = constant$. $F$ is equivalent to a proper polynomial
holomorphic map in $Poly({\BB}^n,{\BB}^{3n-2})$
 if and only if $a = 0$.
\bigskip

\noindent{\it Proof:}\ \ \ \ The projectivization of $F$ is
\[\hat F=\big[tz'(t-\ov a w):twz':w^2\sqrt{1-|a|^2}z': w^2(w-at):
t^2(t-\ov a w)\big].\] Assume $a\neq 0$ and suppose there exist
hyperplanes $H\subset \CC\PP^n$ and $H'\subset \CC\PP^{3n-2}$ such
that $H\cap \ov{{\bf S}^n_1}=\emptyset$, $H'\cap \ov{{\bf
S}^{3n-2}_1} = \emptyset$ and $ \hat F(H) \subset H',\ \ \hat
F(\CC\PP^n \backslash  H) \subset \CC\PP^{3n-2} \backslash H'$.
Then\[\lambda'_1 tz'(t-\ov aw)+\lambda'_2twz' +\lambda'_3
w^2\sqrt{1-|a|^2}z'+\lambda_n w^2(w-at) +\lambda_0 t^2(t-\ov a
w)=(\mu_0t+\mu'z'+\mu_nw)^3\] for some $\lambda'_1, \lambda'_2,
\lambda_3, \mu'\in \CC^{n-1}$ and $\lambda_n,\lambda_0,\mu_0,
\mu_n \in \CC$. \\ Then $\lambda_0=\mu_0^3\neq 0$. We thus can
assume at the beginning that $\lambda_0=\mu=1$.\\ Since there are
no terms like $z_j^3 (j<n) $ in the left hand side, we conclude
that $\mu'=0$. Thus we get
$$\lambda_nw^2(w-at)+t^2(t-\ov a w)=(t+\mu_nw)^3.$$
Therefore $-\ov a=3\mu_n,\ -\lambda_n a=3\mu_n^2,\
\lambda_n=\mu_n^3$ or $\mu_n=-\frac{\ov a}{3}$ and
$\mu_n=-\frac{3}{a}$. This contradicts the assumption that
$0<|a|^2<1$. $\Box$

\end{document}